\documentclass{amsart}
\usepackage{epsfig,psfrag,amsmath,amssymb}
\newtheorem{theorem}{Theorem}

\newtheorem{lemma}{Lemma}

\theoremstyle{definition}
\newtheorem{definition}{Definition}

\theoremstyle{remark}

\newcommand{\begriff}[1]{\textbf{#1}}                  

\newcommand{\llangle}{\langle\langle}
\newcommand{\rrangle}{\rangle\rangle}
\def\co{\colon\thinspace}
\newcommand{\jj}[6]{\left|\begin{smallmatrix} #1 & #2 & #3\\
                                              #6 & #5 & #4
                    \end{smallmatrix}\right|}

\begin{document}

\title{Fast Computation of Secondary Invariants}

\author{Simon A. King}
\address{Simon A. King\\
Mathematisches Forschungsinstitut Oberwolfach\\
Schwarzwaldstr. 9--11\\
D-77709 Oberwolfach\\
Germany}
\email{king@mfo.de}

\begin{abstract}
  A very classical subject in Commutative Algebra is the Invariant Theory of finite groups.
  In our work on $3$--dimensional topology~\cite{KingIdeal2}, we
  found certain examples of group actions on polynomial rings. When we tried to compute the
  invariant ring using \textsc{Singular}~\cite{Singular} or \textsc{Magma}~\cite{magma}, it turned out
  that the existing algorithms did not suffice.
  
  We present here a new algorithm for the computation of secondary invariants, if primary invariants
  are given. Our benchmarks show that the implementation of our algorithm in 
  the library \texttt{finvar} of \textsc{Singular}~\cite{Singular} marks a dramatic improvement in
  the manageable problem size. 
  A particular benefit of our algorithm is that the computation of \emph{irreducible} secondary invariants 
  does not involve the explicit computation of \emph{reducible} secondary invariants, 
  which may save resources.

  The implementation of our algorithm in \textsc{Singular} is for the
  non-modular case; however, the key theorem of our algorithm holds in the modular case as 
  well and might be useful also there.
\end{abstract}

\maketitle

\section{Introduction}

Let $G$ be a finite group, linearly acting on a polynomial ring $R$ with $n$ variables 
over some field $K$. 
We denote the action of $g\in G$ on $r\in R$ by $g.r\in R$. 

Let $R^G=\{r\in R\co g.r=r,\; \forall g\in G\}$ be the invariant ring. Obviously, it is a sub-algebra of $R$,
and one would like to compute generators for $R^G$.
We study here the \begriff{non-modular} case, i.e., the characteristic of $K$ does not divide 
the order of $G$. Note that according to~\cite{Kemper}, algorithms for the non-modular case are 
useful also in the modular case. 

For any subset $S\subset R$, we denote by $\llangle S\rrangle\subset R$ the sub-algebra generated by $S$, 
and by $\langle S\rangle\subset R$ the ideal generated by $S$.
It is well known~\cite{Eisenbud} that there are $n$ (the number of variables) algebraically independent 
homogeneous invariant polynomials $P=\{p_1,...,p_n\}\subset R^G$ such that 
$R^G$ is a finitely generated $\langle\langle P\rangle\rangle$--module.
The elements of $P$ are called \begriff{primary invariants}. Of course, they are not uniquely determined.
There are various algorithms to compute primary invariants~\cite{Kemper}. 
Since the primary invariants are algebraically independent, the sub-algebra $\llangle P\rrangle$
is isomorphic to a polynomial ring with $n$ variables. It is called \begriff{(homogeneous) 
Noetherian normalization} of $R^G$.

Let $S\subset R^G$ be a minimal set of homogeneous $\langle\langle P\rangle\rangle$--module generators of $R^G$. 
The elements of $S$ are called \begriff{secondary invariants}. 
Note that the number of secondary invariants depends on the degrees of the primary invariants. 
Hence, it is advisable to minimize the degrees of the primary invariants.
\begriff{Irreducible} secondary invariants are those non-constant secondary invariants 
that can not be written as a polynomial expression in the primary invariants and the 
other secondary invariants. 
The set of secondary invariants is not unique, even if one fixes the primary invariants. 
It is easy to see that one can choose secondary invariants so that all of them are 
{power products} of irreducible secondary invariants. 

The aim of this paper is to present a new algorithm for the computation of (irreducible) homogeneous 
secondary invariants, if homogeneous primary invariants $P$ are given.
The key theorem for our algorithm concerns \emph{Gr\"obner bases} 
and holds in arbitrary characteristic; however, the algorithm
assumes that we are in the {non-modular} case. 
For simplicity, we even assume that $K$ is of characteristic $0$, but this is not crucial. 

The rest of this paper is organised as follows. In Section~\ref{sec:motivation}, we briefly expose our 
motivating examples arising in low-dimensional topology. 
In the Section~\ref{sec:basic}, we recall the basic scheme
for computing secondary invariants. In Section~\ref{sec:irred_secondary}, we state our key result
and formulate our new algorithm for the computation of (irreducible) secondary invariants. 
In Section~\ref{sec:benchmark}, we provide some examples (partially inspired by our study of problems
in low-dimensional topology) and compare the implementation of
our algorithm in \textsc{Singular}~\cite{Singular} 
with previously implemented algorithms in \textsc{Singular} by 
A. Heydtmann~\cite{Heydtmann} respectively in \textsc{Magma}~\cite{magma} 
by A. Steel~\cite{KemperSteel}.

\section{Motivating examples}
\label{sec:motivation}

The starting point of our work was the study of generalisations of \emph{Turaev--Viro invariants}~\cite{KingIdeal1}, \cite{KingIdeal2}. These are
homeomorphism invariants of compact $3$--dimensional manifolds. Their construction is (with some simplifications) as follows. Let $\mathcal F$ be some finite set, and let $\mathcal T$ be a triangulation of a compact $3$--manifold $M$. An $\mathcal F$--colouring of $\mathcal T$ assigns to any edge of $\mathcal T$ an element of $\mathcal F$. Tetrahedra have six edges. So, for any tetrahedron of $\mathcal T$, an
$\mathcal F$--colouring of $\mathcal T$ gives rise to a six-tuple of colours, that is called \emph{$6j$--symbol}
and denoted by $\jj abcdef$, for $a,b,c,d,e,f\in\mathcal F$. The equivalence classes of $6j$--symbols
with respect to tetrahedral symmetry are variables of some polynomial ring, $R$. The ring also contains
one variable $w_f$ for any $f\in\mathcal F$, called the \emph{weight} of $f$. For any $\mathcal F$--colouring, we
form the product over the weights of the coloured edges and over the $6j$--symbols of the
coloured tetrahedra of $\mathcal T$. By summation over all possible $\mathcal F$--colourings of $\mathcal T$,
we obtain a polynomial $TV(\mathcal T)$ called the \emph{state sum} of $\mathcal T$. Due to the tetrahedral symmetry of the
$6j$--symbols, the state sum is well-defined. However, it depends on the choice of $\mathcal T$ rather than 
on the homeomorphism type of $M$. It was shown by V. Turaev and O. Viro~\cite{TuraevViro} that an appropriate
evaluation of the state sum (yield by the representation theory of Quantum Groups) is independent of the choice of 
$\mathcal T$. This is called a Turaev--Viro invariant. 

In~\cite{KingIdeal1} and~\cite{KingIdeal2}, we define an ideal $I\subset R$, the \emph{Turaev--Viro ideal}.
We show that the coset $tv(M) = TV(\mathcal T) + I$ is independent of $\mathcal T$, hence, a 
homeomorphism invariant of $M$. This generalises the classical Turaev--Viro invariants. By extensive 
computations, we show in~\cite{KingIdeal2} that these so-called \emph{ideal Turaev--Viro invariants} are much stronger than
the classical Turaev--Viro invariants. For this, it was necessary to compute Gr\"obner bases of Turaev--Viro
ideals. It turns out that different algorithms for the computation of Gr\"obner bases differ widely in
their performance. The algorithm \texttt{slimgb} in \textsc{Singular}~\cite{Singular} of M. 
Brickenstein~\cite{Brickenstein} performs particularly well. 

We obtain a lower bound for the number of tetrahedra of any triangulation of $M$, in terms of the minimal degree
of polynomials in the coset $tv(M)$. However, in our computations, the bound appears to be trivial~\cite{KingIdeal2}. There was some hope to improve the lower bound as follows, using computations
of invariant rings.
Let $G$ be the symmetric group of $\mathcal F$. In the obvious way, $G$ acts on the tetrahedral symmetry 
classes of $6j$--symbols and on the weights, and hence, on $R$. The $G$--action permutes the summands
of the state sum. So, the state sum belongs to $R^G$. Let $I^G=I \cap R^G$. Obviously, if $\mathcal T$ is a triangulation
of some compact $3$--manifold $M$, then the coset $tv^G(M)=TV(\mathcal T)+I^G\subset R^G$ is a homeomorphism 
invariant of $M$, and as such in fact equivalent to $tv(M)$. However, since $tv^G(M)\subset tv(M)$, there is some hope that the minimal
degree of polynomials in the coset $tv^G(M)$ is higher than in $tv(M)$, which would provide stronger bounds
for the number of tetrahedra.

This is how we became interested in the computation of invariant rings.
The existing implementations in \textsc{Magma} and \textsc{Singular} could not
compute the secondary invariants in several of our examples. This motivated us to develop a new algorithm
for the computation of secondary invariants. It has been part of the \texttt{finvar} library
of \textsc{Singular}~\cite{Singular} since release 3-0-2 (July 2006).
Unfortunately, in our topological applications, we did not find an improvement of the lower bound for the number of tetrahedra. However, our new algorithm for the computation of secondary invariants certainly is  of independent
interest.

\section{Generalities on the computation of secondary invariants}
\label{sec:basic}

In the non-modular case, we can use the Reynolds operator $\mathrm{Rey}\co R\to R^G$,
which is defined by 
\[ \mathrm{Rey}(r) = \frac 1{|G|}\sum_{g\in G} g.r \]
for $r\in R$.
By construction, the restriction of the Reynolds operator to $R^G$ is the identity.
Let $B_d\subset R^G$ be the images under the Reynolds
operator of all monomials of $R$ of degree $d$. 
It is well known that one can find a system of homogeneous secondary invariants of degree 
$d$ in $B_d$~\cite{Sturmfels}. 
But how can one determine what elements of $B_d$ are eligible as secondary invariants?

Let $S_0,S_1,S_2,...,S_{d-1}\subset R^G$ be the homogeneous secondary invariants of degree $0,1,2,...,d-1$,
respectively (we can take $S_0=\{1\}$), and let $IS_i\subset S_i$ be the irreducible ones, for $i=1,...,d-1$. 
Let $s_1,...,s_m\in R^G$ be some homogeneous secondary invariants of degree $d$. Let $b\in B_d$. We can 
choose $b$ as a new homogeneous secondary invariant of $R^G$, if $b$ is not contained
in the $\llangle P\rrangle$--module generated by $S_0\cup S_1\cup \cdots \cup S_{d-1}\cup \{s_1,...,s_m\}$.
It is not difficult to show that this is the case if and only if
$b$ is not contained in the \emph{ideal} $\langle P\cup \{s_1,...,s_m\}\rangle\subset R$; see~\cite{Sturmfels}.

Ideal membership can be tested using \emph{Gr\"obner bases}. 
For $p\in R$ and a finite subset $\mathcal G\subset R$, we denote the \begriff{remainder} of $p$ by reduction
modulo $\mathcal G$ by $\mathrm{rem}(b;\mathcal G)$. The remainder is iteratively defined,
depends on the choice of a monomial order, and in general depends on the order of the elements 
of $\mathcal G$. For a definition of remainder, of Gr\"obner bases, 
and for a proof of the following classical result, we refer to~\cite{Froberg} or~\cite{KreuzerRobbiano}.
\begin{theorem}\label{thm:idealmembership}
  Let $\mathcal G$ be a Gr\"obner basis of $\langle \mathcal G\rangle\subset R$, and let $p\in R$.
  Then, $\mathrm{rem}(p;\mathcal G)$ does not depend on the order of polynomials in $\mathcal G$,
  and we have $\mathrm{rem}(p;\mathcal G)=0$ if and only if $p\in\langle \mathcal G\rangle$.\qed
\end{theorem}

We thus obtain the following very basic algorithm for finding homogeneous secondary invariants $S_d$ 
of degree $d$, provided those of smaller degrees have been computed before.
\\
\textsc{Basic Algorithm}
\begin{enumerate}
\item Let $S_d=\emptyset$. Let $\mathcal G$ be a Gr\"obner basis of $\langle P\rangle$.
\item For all $b\in B_d$:\\
If $b\not\in \langle P\cup S_d\rangle$ (which is tested by reduction modulo $\mathcal G$)
 then replace $S_d$ by $S_d\cup \{b\}$; compute a Gr\"obner basis
 of $\langle P\cup S_d\rangle$ and replace $\mathcal G$ with it.
\item Return $S_d$.
\end{enumerate}

There are several ways to improve this algorithm. One way is an application of Molien's 
Theorem~\cite{Sturmfels}, \cite{Kemper}, \cite{Heydtmann}. We will not go into details here. 
Molien's Theorem allows to 
compute the number $m_d$ of secondary invariants of degree $d$. In other words, if in the 
above algorithm we got $m_d$ secondary invariants, we can immediately 
break the loop in Step~(2).

We also would  like to see which of the secondary invariants in $S_d$ are irreducible, since these, together
with $P$, generate $R^G$ as a sub-algebra of $R$. For that 
purpose, one forms all power products of degree $d$ of elements of $IS_1\cup IS_2\cup \cdots \cup IS_{d-1}$
and chooses from them as many secondary invariants as possible
(compare~\cite{Heydtmann} or~\cite{Kemper}). 
If there are further secondary invariants (which we know from computation of $m_d$),
then one proceeds as above with $B_d$, and obtains all \emph{irreducible} secondary invariants $IS_d$ of
degree $d$. So, the algorithm is as follows.
\\
\textsc{Refined Algorithm}
\begin{enumerate}
\item Compute $m_d$. Let $S_d=IS_d=\emptyset$ and let $\mathcal G$ be a Gr\"obner basis of 
      $\langle P\rangle$.
\item For all power products $b$ of degree $d$ of elements of $IS_1\cup IS_2\cup \cdots \cup IS_{d-1}$:
   \begin{enumerate}
   \item If $b\not \in \langle P\cup S_d\rangle$ (which is tested using $\mathcal G$)
      then replace $S_d$ by $S_d\cup \{b\}$; compute a Gr\"obner basis
      of $\langle P\cup S_d\rangle$ and replace $\mathcal G$ with it.
   \item If $|S_d|=m_d$ then break and return $(S_d,IS_d)$.
   \end{enumerate}
\item For all $b\in B_d$:
   \begin{enumerate}
   \item If $b\not \in \langle P\cup S_d\rangle$ (which is tested using $\mathcal G$)
      then replace $S_d$ by $S_d\cup \{b\}$, and $IS_d$ by $IS_d\cup \{b\}$; compute a Gr\"obner basis
      of $\langle P\cup S_d\rangle$ and replace $\mathcal G$ with it.
   \item If $|S_d|=m_d$ then break and return $(S_d,IS_d)$.
   \end{enumerate}
\end{enumerate}

Eventually, $S_d$ contains homogeneous secondary invariants of degree $d$, and $IS_d$ contains the irreducible
ones.
In this form, the algorithm has been implemented in 1998 by A. Heydtmann~\cite{Heydtmann} as the 
procedure \texttt{secondary\underline{\ }char0} of the library \texttt{finvar} of \textsc{Singular}. 
In Step~(2), 
the ideal membership is tested by computing the remainder modulo some Gr\"obner basis of the ideal.
This ideal changes once a new secondary invariant has been found. So, the algorithm involves 
many Gr\"obner basis computations. This is its main disadvantage and
limits the applicability of the Basic and the Refined Algorithm.

An alternative algorithm was proposed by Kemper and Steel
(see~\cite{Kemper},~\cite{KemperSteel} or~\cite{DerksenKemper})
and implemented in \textsc{Magma}~\cite{magma}. 
Here, new secondary invariants are detected not by a general solution of the ideal membership problem but
by Linear Algebra. 
This algorithm only involves one Gr\"obner basis computation, namely for the ideal $\langle P\rangle$.
But for computing some of the invariant rings that arise 
in our study of homeomorphism invariants of $3$--dimensional manifolds~\cite{KingIdeal2}, 
this does not suffice either.

\section{The New Algorithm}
\label{sec:irred_secondary}

The main feature of our new algorithm is that, after computing some (homogeneous) Gr\"obner 
basis of $\langle P\rangle$, 
we can directly write down a \emph{homogeneous Gr\"obner basis
up to degree $d$} of $\langle P\cup S_d\rangle$, once a new secondary invariant of degree $d$ has been found. 
We can do so whithout any lengthy computations
(in contrast to~\cite{Heydtmann}), and we also avoid to deal with 
huge systems of linear equations (in contrast to~\cite{KemperSteel}, \cite{Kemper},
\cite{DerksenKemper}). 
This allows to solve the ideal membership problem in a very quick way. 
We recall the notion of \lq\lq homogeneous Gr\"obner bases up to 
degree $d$\rq\rq\ in the following paragraphs. At the end of the section, we 
provide our key theorem and formulate our new algorithm.

For $p\in R$, let $lm(p)$ the \begriff{leading monomial} of $p$, let $lc(p)$ be the coefficient 
of $lm(p)$ in $p$, and let $lt(p)=lc(p)lm(p)$ be the \begriff{leading term} of $p$. 
The least common multiple is denoted by $LCM(\cdot,\cdot)$.
Now we can recall the definition of the \begriff{$S$--polynomial} of $p,q\in R$:
\[S(p,q) = \frac{LCM(lm(p),lm(q))}{lt(p)} p - \frac{LCM(lm(p),lm(q))}{lt(q)} q\]
Obviously, the $S$--polynomial of $p$ and $q$ belongs to the ideal $\langle p,q\rangle\subset R$.
The leading terms of $p$ and $q$ are canceling one another, so, the leading monomial
of $S(p,q)$ corresponds to monomials of $p$ or $q$ that are not leading. 
The following result can be found, e.g., in~\cite{Froberg} or~\cite{KreuzerRobbiano}.
\begin{theorem}[Buchberger's Criterion]\label{thm:buchberger}
  A set $g_1,...,g_k\in R$ of polynomials is a Gr\"obner basis of the 
  ideal $\langle g_1,...,g_k\rangle\subset R$ if and only if
  $\mathrm{rem}\left(S(g_i,g_j); g_1,...,g_k\right) = 0$ for all $i,j=1,...,k$.\qed
\end{theorem}

Buchberger's Criterion directly leads to Buchberger's algorithm for the construction of
a Gr\"obner basis of an ideal: One starts with any generating set of the ideal. If the remainder
modulo the generators of the $S$--polynomial
of some pair of generators does not vanish, then the remainder 
is added as a new generator. 
This will be repeated until all $S$--polynomials reduce to $0$; it can be shown that this will eventually be 
the case, after finitely many steps.

Here, we are in a special situation: We work with homogeneous polynomials. It is easy to see that
if $p$ and $q$ are homogeneous then so is $S(p,q)$, and its degree is higher than the maximum of the
degrees of $p$ and $q$, unless $lm(p)=lm(q)$. 
If $p, g_1,g_2,...,g_k\in R$ are homogeneous then so is $\mathrm{rem}(p;g_1,...,g_k)$.
Moreover, either $\mathrm{rem}(p;g_1,...,g_k)= 0$ or $\deg\left(\mathrm{rem}(p;g_1,...,g_k)\right)
= \deg (p)$. 
For computing $\mathrm{rem}(p;g_1,...,g_k)$, only those $g_i$ play a role with 
$\deg(g_i)\le \deg(p)$, for $i=1,...,k$.
It follows: If an ideal $I\subset R$ is homogeneous (i.e., it can be generated by homogeneous polynomials) 
then it has a Gr\"obner basis of homogeneous polynomials. Such a Gr\"obner basis can be 
constructed degree-wise.

\begin{definition}
  A finite set $\{g_1,...,g_k\}\subset R$ of homogeneous polynomials is a \begriff{homogeneous 
  Gr\"obner basis up to degree $d$}
  of the ideal $\langle g_1,...,g_k\rangle$, if 
  $$\mathrm{rem}\left(S(g_i,g_j); g_1,...,g_k\right) = 0$$ or $\deg\left(S(g_i,g_j)\right)>d$,
  for all $i,j=1,...,k$.
\end{definition}

\begin{lemma}\label{lem:idealmembership}
  Let $\{g_1,...,g_k\}\subset R$ be a homogeneous Gr\"obner basis up to degree $d$, and let $p\in R$ be
  a homogeneous polynomial of degree at most $d$. Then, $p\in\langle g_1,...,g_k\rangle$ if and 
  only if $\mathrm{rem}\left(p;g_1,...,g_k\right) =0$.
\end{lemma}
\begin{proof}
  The paragraph preceding the definition implies that $\{g_1,...,g_k\}$ can be extended to 
  a Gr\"obner basis $\mathcal G$ of $\langle g_1,...,g_k\rangle$ by adding homogeneous polynomials whose 
  degrees exceed $d$. Since $\deg(p)\le d$, we have $\mathrm{rem}(p;\mathcal G)=\mathrm{rem}(p;g_1,...,g_k)$.
  Since $p\in\langle \mathcal G\rangle$ if and only if $\mathrm{rem}(p;\mathcal G)=0$ by 
  Theorem~\ref{thm:idealmembership}, the result follows.
\end{proof}

We see that in order to do Step~(2) in the \textsc{Basic Algorithm} (or the corresponding steps in the 
\textsc{Refined Algorithm}) it suffices to know a homogeneous Gr\"obner basis 
up to degree $d$ of $\langle P\cup S_d\rangle$.
Our key theorem states that this Gr\"obner basis can be constructed iteratively, as follows.
\begin{theorem}\label{thm:keythm}
  Let $\mathcal G\subset R$ be a homogeneous Gr\"obner basis up to degree $d$ of $\langle \mathcal G\rangle$.
  Let $p\in R$ be a homogeneous polynomial of degree $d$, and $p\not\in \langle \mathcal G\rangle$.
  Then $\mathcal G\cup \{\mathrm{rem}(p;\mathcal G)\}$ is a homogeneous Gr\"obner basis up 
  to degree $d$ of $\langle \mathcal G\cup\{p\}\rangle$.
\end{theorem}
\begin{proof}
  Let $r=\mathrm{rem}(p;\mathcal G)$. Since $p\not\in \langle \mathcal G\rangle$ and all polynomials
  are homogeneous, we 
  have $r\not=0$, $\deg(r)=d$, and
  $\langle \mathcal G\cup\{p\}\rangle = \langle \mathcal G\cup
  \{r\}\rangle$.

  By hypothesis, the $S$--polynomials of pairs of elements of $\mathcal G$ are of degree $>d$ or reduce to
  $0$ modulo $\mathcal G$. 
  We now consider the $S$--polynomials of $r$ and elements of $\mathcal G$. 
  Let $g\in \mathcal G$.  
  By definition of the remainder, we have $lm(r) \not= lm(g)$. 
  Therefore the $S$--polynomial of $r$ and $g$ is of degree $>d=\deg(r)$. Thus the claim
  follows.
\end{proof}

\noindent
We obtain the 
\\
\textsc{New Algorithm}
\begin{enumerate}
\item Compute $m_d$ and a homogeneous Gr\"obner basis $\mathcal G$ of 
   $\langle P\rangle$. Let $S_d=IS_d=\emptyset$.
\item For all power products $b$ of degree $d$ of elements of $IS_1\cup IS_2\cup \cdots \cup IS_{d-1}$:
   \begin{enumerate}
   \item If $\mathrm{rem}(b;\mathcal G)<>0$
      then replace $S_d$ by $S_d\cup \{b\}$ and $\mathcal G$ by $\mathcal G\cup \{\mathrm{rem}(b;\mathcal G)\}$.
   \item If $|S_d|=m_d$ then break and return $(S_d,IS_d)$.
   \end{enumerate}
\item For all $b\in B_d$:
   \begin{enumerate}
   \item If $\mathrm{rem}(b;\mathcal G)<>0$
      then replace $S_d$ by $S_d\cup \{b\}$, $IS_d$ by $IS_d\cup \{b\}$
      and $\mathcal G$ by $\mathcal G\cup \{\mathrm{rem}(b;\mathcal G)\}$.
   \item If $|S_d|=m_d$ then break and return $(S_d,IS_d)$.
   \end{enumerate}
\end{enumerate}
By Theorem~\ref{thm:keythm} and induction, $\mathcal G$ is a homogeneous Gr\"obner basis up to degree $d$
of $\langle P\cup S_d\rangle$. Hence, in Step~(2)(a) and~(3)(a) one has $\mathrm{rem}(b;\mathcal G)<>0$ if
and only if $b\not\in \langle P\cup S_d\rangle$.
The \textsc{New Algorithm} is a dramatic improvement of the \textsc{Refined Algorithm}.
However, in our examples this was still not enough.

One should take more care in Step~(2) of the \textsc{New Algorithm}. It simply says 
\lq\lq For all power products $b$ of degree $d$ of elements of $IS_1\cup IS_2\cup \cdots \cup IS_{d-1}$\rq\rq.
Two questions arise: 
\begin{enumerate}
\item How shall one generate the power products? 
\item Is it necessary to generate \emph{all} possible power products, or can one restrict the search?
\end{enumerate}

In very complex computations, the number of power products is gigantic. But usually only a small proportion of 
them will be eligible as secondary invariant.
So, for saving computer's memory, it is advisable to generate the power products one after the other (or in small packages), rather than generating all power products at once; this answers Question~(1).

Apparently Question~(2) was never addressed in the literature. However, it turns
out that a careful choice of power products provides another dramatic improvement
of the performance of the algorithm.
Our choice is based on the following lemma. This lemma seems to be well known, but
to the best of the author's knowledge it did not appear in the literature and it
was not used in implementations.

\begin{lemma}\label{lem:RestrictChoice}
  Assume that secondary invariants of degree $<d$ are computed such that all of them are 
  power products of irreducible secondary invariants.
  In the quest for reducible homogeneous secondary invariants of degree $d$, it suffices to consider
  power products of the form $i\cdot s$, where $i$ is a homogeneous irreducible secondary invariant of degree
  $<d$, and $s$ is some secondary invariant of degree $d-\deg(i)$.
\end{lemma}
\begin{proof}
  Let $p\in R$ be a power product of degree $d$ of irreducible secondary invariants. Hence, it can be written
  as $p=iq$, with an irreducible homogeneous secondary invariant $i$ of degree $<d$ and some homogeneous
  $G$--invariant polynomial $q$ of degree $d-\deg (i)$ (we do not use that $q$ is a power product of irreducible
  secondary invariants). 

  Recall that the secondary invariants generate the invariant ring as a $\llangle P\rrangle$--module. 
  Hence one can rewrite
  $q = q_0 + k_1s_1+\cdots + k_ts_t$, where $q_0\in\langle P\rangle$,
  $k_1,...,k_t\in K$, and $s_1,...s_t$ are homogeneous secondary invariants of degree $\deg(q)$.
  We obtain $p=iq_0 + k_1(is_1)+\cdots + k_t(is_t)$. Hence, rather than chosing $p$ as a 
  $\llangle P\rrangle$--module generator of $R^G$, we may choose $is_1,...,is_t$, which, by induction, are
  all power products of irreducible secondary invariants.
\end{proof}

\noindent
\textsc{Improved New Algorithm}
\begin{enumerate}
\item Compute $m_d$. Let $\mathcal G$ be a Gr\"obner basis of $\langle P\rangle$. Let $S_d=IS_d=\emptyset$.
\item For all products $b=i\cdot s$ with $i\in IS_1\cup\cdots IS_{d-1}$ and $s\in S_{d-\deg(i)}$:
   \begin{enumerate}
   \item If $\mathrm{rem}(b;\mathcal G)\not=0$
      then replace $S_d$ by $S_d\cup \{b\}$ and $\mathcal G$ by $\mathcal G\cup \{\mathrm{rem}(b;\mathcal G)\}$.
   \item If $|S_d|=m_d$ then break and return $(S_d,IS_d)$.
   \end{enumerate}
\item For all $b\in B_d$:
   \begin{enumerate}
   \item If $\mathrm{rem}(b;\mathcal G)\not=0$
      then replace $S_d$ by $S_d\cup \{b\}$, $IS_d$ by $IS_d\cup \{b\}$
      and $\mathcal G$ by $\mathcal G\cup \{\mathrm{rem}(b;\mathcal G)\}$.
   \item If $|S_d|=m_d$ then break and return $(S_d,IS_d)$.
   \end{enumerate}
\end{enumerate}

This is the algorithm that is implemented as \texttt{secondary\underline{\ }char0} in the 
library \texttt{finvar} of \textsc{Singular}
3-0-2~\cite{Singular}, released in Juli 2006. In Step~(2), the secondary invariant
$s$ may be a non-trivial powerproduct itself, hence, can be expressed as $s=i_ss'$, where $i_s$ is an irreducible 
secondary invariant and $s'$ is (by induction) some other secondary invariant. 
Of course one should consider only one of the two products $i_s(is')$ and $i(i_ss')$ in the enumeration. 

Often one is only interested in the irreducible secondary invariants, which, together with the primary 
invariants, generate the invariant ring as a sub-algebra. Therefore we 
implemented yet another version of 
the \textsc{Improved New Algorithm} in \textsc{Singular} 3-0-2, 
namely \texttt{irred\underline{\ }secondary\underline{\ }char0}. This algorithm computes irreducible 
secondary invariants, but does not explicitely compute the reducible secondary invariants. 
That works as follows.

Let $\mathcal G_P$ be a Gr\"obner basis of $\langle P\rangle$.  
In Step~(2)(a) of the \textsc{Improved New Algorithm}, one replaces $S_d$ by $S_d\cup \{\mathrm{rem}(b;\mathcal G_P)\}$, rather than by 
$S_d\cup \{b\}$. In Step~(3)(a) one replaces $S_d$ by $S_d\cup \{\mathrm{rem}(b;\mathcal G_P)\}$ and $IS_d$ by $IS_d\cup\{b\}$.
In the end, $S_d$ does not contain secondary invariants, but 
\emph{normal forms} of secondary invariants with respect to $\mathcal G_P$.
Since $\mathrm{rem}\left(\mathrm{rem}(p_1;\mathcal G_P)\cdot\mathrm{rem}(p_2;\mathcal G_P);
\mathcal G_P\right) = \mathrm{rem}(p_1\cdot p_2;\mathcal G_P)$ and since a reduction modulo $\mathcal G$
in Steps~(2)(a) and~(3)(a) also comprises a reduction modulo $\mathcal G_P$,
this maintains all informations that one needs for determining how many secondary invariants are reducible
in Step~(2) and for finding the irreducible secondary invariants in Step~(3). So in the end, $IS_d$
contains the irreducible secondary invariants in degree $d$.
This detail of our implementation very often saves much memory and computation time, as can be seen 
in Table~\ref{tab:bench} in Examples (1) and (6)--(9).
In Example~(8), we can compute the irreducible secondary invariants although 
the computation of all 31104 secondary invariants exceeds the resources.

An example of Kemper (example~(9) in the next Section) motivated us to further refine the implementation of the \textsc{Improved New Algorithm}. It concerns the generation of $B_d$: If there are irreducible secondary invariants in rather high degrees $d$ (in Kemper's example, there are two irreducible secondary invariants of degree $9$), it is advisable to generate not all of $B_d$ at once, but in small portions. This will be part of release 3-0-3 of \textsc{Singular}.

\section{Benchmark Tests for the Computation of Invariant Rings}
\label{sec:benchmark}

\subsection{The Test Examples}
We already mentioned that some of our test examples arise in low-dimensional topology. This yields
Examples~(1), (7) and (8). For background information, see~\cite{KingIdeal2}.
We will not go into details here, but 
just provide the matrices and primary invariants of our nine test examples.
They are roughly ordered by increasing computation time. 
The ring variables are called $x_1,x_2,...$. Let $e_i$ be the column vector with $1$ in
position $i$ and $0$ otherwise. 
Our focus was not on the computation of primary invariants; note that in various examples
the primary invariants are not optimal. 
\begin{enumerate}
\item A $13$--dimensional representation of the symmetric group $S_2$ is given by the matrix
   \[
    M = \left(e_{2}e_{1}e_{13}e_{12}e_{11}e_{8}e_{10}e_{6}e_{9}e_{7}e_{5}e_{4}e_{3}\right)
%    \begin{smallmatrix}
%0& 1& 0& 0& 0& 0& 0& 0& 0& 0& 0& 0& 0\\
%1& 0& 0& 0& 0& 0& 0& 0& 0& 0& 0& 0& 0\\
%0& 0& 0& 0& 0& 0& 0& 0& 0& 0& 0& 0& 1\\
%0& 0& 0& 0& 0& 0& 0& 0& 0& 0& 0& 1& 0\\
%0& 0& 0& 0& 0& 0& 0& 0& 0& 0& 1& 0& 0\\
%0& 0& 0& 0& 0& 0& 0& 1& 0& 0& 0& 0& 0\\
%0& 0& 0& 0& 0& 0& 0& 0& 0& 1& 0& 0& 0\\
%0& 0& 0& 0& 0& 1& 0& 0& 0& 0& 0& 0& 0\\
%0& 0& 0& 0& 0& 0& 0& 0& 1& 0& 0& 0& 0\\
%0& 0& 0& 0& 0& 0& 1& 0& 0& 0& 0& 0& 0\\
%0& 0& 0& 0& 1& 0& 0& 0& 0& 0& 0& 0& 0\\
%0& 0& 0& 1& 0& 0& 0& 0& 0& 0& 0& 0& 0\\
%0& 0& 1& 0& 0& 0& 0& 0& 0& 0& 0& 0& 0
%              \end{smallmatrix}\right)
   \]
   Our primary invariants are
   \begin{eqnarray*}
     &&x_{9},\;\;x_{7}+x_{10},\;\;x_{6}+x_{8},\;\;x_{5}+x_{11},\;\;
     x_{4}+x_{12},\;\;x_{3}+x_{13},\\ &&x_1+x_2,\;\;x_{3}x_{13},\;\;
     x_{4}x_{12},\;\;x_{5}x_{11},\;\;x_{7}x_{10},\;\;x_{6}x_{8},\;\;x_1x_2
   \end{eqnarray*}
   There are 32 secondary invariants of maximal degree $6$, among which are $15$ irreducible 
   secondary invariants up to degree $2$.
\item A $6$--dimensional representation of $S_4$ is given by the matrices
   \begin{eqnarray*}
    M_1 &=& \left(e_{1}e_{4}e_{5}e_{2}e_{3}e_{6}\right)
%    \begin{smallmatrix}
%      1& 0& 0& 0& 0& 0\\
%      0& 0& 0& 1& 0& 0\\
%      0& 0& 0& 0& 1& 0\\
%      0& 1& 0& 0& 0& 0\\
%      0& 0& 1& 0& 0& 0\\
%      0& 0& 0& 0& 0& 1
%              \end{smallmatrix}\right) 
              \\
    M_2 &=& \left(e_{4}e_{1}e_{5}e_{2}e_{6}e_{3}\right)
%    \begin{smallmatrix}
%      0& 1& 0& 0& 0& 0\\
%      0& 0& 0& 1& 0& 0\\
%      0& 0& 0& 0& 0& 1\\
%      1& 0& 0& 0& 0& 0\\
%      0& 0& 1& 0& 0& 0\\
%      0& 0& 0& 0& 1& 0
%              \end{smallmatrix}\right)
   \end{eqnarray*}
   Our primary invariants are
   \begin{eqnarray*}
    &&x_3+x_5+x_6,\;\;x_1+x_2+x_4,\;\;x_3x_5+x_3x_6+x_5x_6,\\
    &&x_3x_4+x_2x_5+x_1x_6,\;\; x_1x_2x_4,\;\;
    x_1^3x_2^3+x_1^3x_4^3+x_2^3x_4^3+x_3^2x_5^2x_6^2
   \end{eqnarray*}
   There are $12$ secondary invariants of maximal degree $9$, among which are $4$ irreducible
   secondary invariants of maximal degree $3$.

\item A $6$--dimensional representation of the alternating group $A_4$ is given by the matrices
   \begin{eqnarray*}
    M_1 &=& \left(e_{4}e_{1}e_{5}e_{2}e_{6}e_{3}\right)
%    \begin{smallmatrix}
%  0& 1& 0& 0& 0& 0\\
%  0& 0& 0& 1& 0& 0\\
%  0& 0& 0& 0& 0& 1\\
%  1& 0& 0& 0& 0& 0\\
%  0& 0& 1& 0& 0& 0\\
%  0& 0& 0& 0& 1& 0
%              \end{smallmatrix}\right) 
              \\
    M_2 &=& \left(e_{2}e_{3}e_{1}e_{6}e_{4}e_{5}\right)
%    \begin{smallmatrix}
%  0& 0& 1& 0& 0& 0\\
%  1& 0& 0& 0& 0& 0\\
%  0& 1& 0& 0& 0& 0\\
%  0& 0& 0& 0& 1& 0\\
%  0& 0& 0& 0& 0& 1\\
%  0& 0& 0& 1& 0& 0
%              \end{smallmatrix}\right)
   \end{eqnarray*}
   Our primary invariants are
   \begin{eqnarray*}
     && x_1+x_2+x_3+x_4+x_5+x_6,\;\;
      x_3x_4+x_2x_5+x_1x_6,\\
     && x_1x_2+x_1x_3+x_2x_3+x_1x_4+x_2x_4+x_1x_5+x_3x_5\\
     && \qquad\mbox{} +x_4x_5+x_2x_6+x_3x_6+x_4x_6+x_5x_6,\\
     && x_3^2x_4+x_3x_4^2+x_2^2x_5+x_2x_5^2+x_1^2x_6+x_1x_6^2,\\
     && x_1x_2x_4+x_1x_3x_5+x_2x_3x_6+x_4x_5x_6, \\
     && x_1^2x_2^4+x_1^4x_3^2+x_2^2x_3^4+x_1^4x_4^2+x_2^2x_4^4\\
     && \qquad\mbox{}+x_3^4x_5^2+x_4^4x_5^2+x_1^2x_5^4+x_2^4x_6^2+x_5^4x_6^2+x_3^2x_6^4+x_4^2x_6^4
   \end{eqnarray*}
   There are 18 secondary invariants of maximal degree 11, among which are 8 irreducible
   secondary invariants of maximal degree 5.

\item A $6$--dimensional representation of the dihedral group $D_6$ is given by the matrices
   \begin{eqnarray*}
    M_1 &=& \left(e_{6}e_{5}e_{4}e_{3}e_{2}e_{1}\right)
%    \begin{smallmatrix}
%0&0&0&0&0&1\\
%0&0&0&0&1&0\\
%0&0&0&1&0&0\\
%0&0&1&0&0&0\\
%0&1&0&0&0&0\\
%1&0&0&0&0&0
%              \end{smallmatrix}\right) 
              \\
    M_2 &=& \left(e_{3}e_{1}e_{2}e_{6}e_{4}e_{5}\right)
%    \begin{smallmatrix}
%0&1&0&0&0&0\\
%0&0&1&0&0&0\\
%1&0&0&0&0&0\\
%0&0&0&0&1&0\\
%0&0&0&0&0&1\\
%0&0&0&1&0&0
%              \end{smallmatrix}\right) 
   \end{eqnarray*}
   Our primary invariants are the elementary symmetric polynomials.
   There are 120 secondary invariants of maximal degree 14, among which are 10 irreducible
   secondary invariants of maximal degree 4.

\item A $8$--dimensional representation of $D_8$ is given by the matrices
   \begin{eqnarray*}
    M_1 &=& \left(e_{8}e_{7}e_{6}e_{5}e_{4}e_{3}e_{2}e_{1}\right)
%    \begin{smallmatrix}
%0& 0& 0& 0& 0& 0& 0& 1\\
%0& 0& 0& 0& 0& 0& 1& 0\\
%0& 0& 0& 0& 0& 1& 0& 0\\
%0& 0& 0& 0& 1& 0& 0& 0\\
%0& 0& 0& 1& 0& 0& 0& 0\\
%0& 0& 1& 0& 0& 0& 0& 0\\
%0& 1& 0& 0& 0& 0& 0& 0\\
%1& 0& 0& 0& 0& 0& 0& 0
%              \end{smallmatrix}\right) 
              \\
    M_2 &=& \left(e_{4}e_{1}e_{2}e_{3}e_{8}e_{5}e_{6}e_{7}\right)
%    \begin{smallmatrix}
%0& 1& 0& 0&  0& 0& 0& 0\\
%0& 0& 1& 0&  0& 0& 0& 0\\
%0& 0& 0& 1&  0& 0& 0& 0\\
%1& 0& 0& 0&  0& 0& 0& 0\\
%0& 0& 0& 0&  0& 1& 0& 0\\
%0& 0& 0& 0&  0& 0& 1& 0\\
%0& 0& 0& 0&  0& 0& 0& 1\\
%0& 0& 0& 0&  1& 0& 0& 0\\
%              \end{smallmatrix}\right) 
   \end{eqnarray*}
   Our primary invariants are
   \begin{eqnarray*}
&&x_1+x_2+x_3+x_4+x_5+x_6+x_7+x_8,\\
&&x_4x_5+x_1x_6+x_2x_7+x_3x_8,\;\;
x_3x_5+x_4x_6+x_1x_7+x_2x_8,\\
&&x_2x_5+x_3x_6+x_4x_7+x_1x_8,\;\;
x_1x_5+x_2x_6+x_3x_7+x_4x_8,\\
&&x_1x_3+x_2x_4+x_5x_7+x_6x_8,\;\;
x_1x_2x_3x_4+x_5x_6x_7x_8,\\
&&x_1x_2^3+x_2x_3^3+x_1^3x_4+x_3x_4^3
%&&\qquad\mbox{}
+x_5^3x_6+x_6^3x_7+x_7^3x_8+x_5x_8^3
   \end{eqnarray*}
   There are 64 secondary invariants of maximal degree 11, among which are 24 irreducible
   secondary invariants of maximal degree 5.

\item A $7$--dimensional representation of $D_{14}$ is given by the matrices
   \begin{eqnarray*}
    M_1 &=& \left(e_{2}e_{3}e_{4}e_{5}e_{6}e_{7}e_{1}\right)
%    \begin{smallmatrix}
%0&0&0&0&0&0&1\\
%1&0&0&0&0&0&0\\
%0&1&0&0&0&0&0\\
%0&0&1&0&0&0&0\\
%0&0&0&1&0&0&0\\
%0&0&0&0&1&0&0\\
%0&0&0&0&0&1&0 
%              \end{smallmatrix}\right) 
              \\
    M_2 &=& \left(e_{1}e_{7}e_{6}e_{5}e_{4}e_{3}e_{2}\right)
%    \begin{smallmatrix}
%1&0&0&0&0&0&0\\
%0&0&0&0&0&0&1\\
%0&0&0&0&0&1&0\\
%0&0&0&0&1&0&0\\
%0&0&0&1&0&0&0\\
%0&0&1&0&0&0&0\\
%0&1&0&0&0&0&0
%              \end{smallmatrix}\right) 
   \end{eqnarray*}
   Our primary invariants are  the elementary symmetric polynomials.
   There are 360 secondary invariants of maximal degree 18, among which are 19 irreducible
   secondary invariants of maximal degree 7.

\item A $15$--dimensional representation of $S_3$ is given by the matrices
   \begin{eqnarray*}
    M_1 &=& \left(e_{2}e_{1}e_{3}e_{4}e_{7}e_{14}e_{5}e_{8}e_{11}e_{13}e_{9}e_{15}e_{10}e_{6}e_{12}\right)
%    \begin{smallmatrix}
%0&1&0&0&0&0&0&0&0&0&0&0&0&0&0\\
%1&0&0&0&0&0&0&0&0&0&0&0&0&0&0\\
%0&0&1&0&0&0&0&0&0&0&0&0&0&0&0\\
%0&0&0&1&0&0&0&0&0&0&0&0&0&0&0\\
%0&0&0&0&0&0&1&0&0&0&0&0&0&0&0\\
%0&0&0&0&0&0&0&0&0&0&0&0&0&1&0\\
%0&0&0&0&1&0&0&0&0&0&0&0&0&0&0\\
%0&0&0&0&0&0&0&1&0&0&0&0&0&0&0\\
%0&0&0&0&0&0&0&0&0&0&1&0&0&0&0\\
%0&0&0&0&0&0&0&0&0&0&0&0&1&0&0\\
%0&0&0&0&0&0&0&0&1&0&0&0&0&0&0\\
%0&0&0&0&0&0&0&0&0&0&0&0&0&0&1\\
%0&0&0&0&0&0&0&0&0&1&0&0&0&0&0\\
%0&0&0&0&0&1&0&0&0&0&0&0&0&0&0\\
%0&0&0&0&0&0&0&0&0&0&0&1&0&0&0
%              \end{smallmatrix}\right) 
              \\
    M_2 &=& \left(e_{1}e_{3}e_{2}e_{4}e_{5}e_{9}e_{8}e_{7}e_{6}e_{13}e_{12}e_{11}e_{10}e_{15}e_{14}\right)
%    \begin{smallmatrix}
%1&0&0&0&0&0&0&0&0&0&0&0&0&0&0\\
%0&0&1&0&0&0&0&0&0&0&0&0&0&0&0\\
%0&1&0&0&0&0&0&0&0&0&0&0&0&0&0\\
%0&0&0&1&0&0&0&0&0&0&0&0&0&0&0\\
%0&0&0&0&1&0&0&0&0&0&0&0&0&0&0\\
%0&0&0&0&0&0&0&0&1&0&0&0&0&0&0\\
%0&0&0&0&0&0&0&1&0&0&0&0&0&0&0\\
%0&0&0&0&0&0&1&0&0&0&0&0&0&0&0\\
%0&0&0&0&0&1&0&0&0&0&0&0&0&0&0\\
%0&0&0&0&0&0&0&0&0&0&0&0&1&0&0\\
%0&0&0&0&0&0&0&0&0&0&0&1&0&0&0\\
%0&0&0&0&0&0&0&0&0&0&1&0&0&0&0\\
%0&0&0&0&0&0&0&0&0&1&0&0&0&0&0\\
%0&0&0&0&0&0&0&0&0&0&0&0&0&0&1\\
%0&0&0&0&0&0&0&0&0&0&0&0&0&1&0
%              \end{smallmatrix}\right) 
   \end{eqnarray*}
   Our primary invariants are
   \begin{eqnarray*}
&& x_{1}+x_{2}+x_{3},\;\;x_{1}x_{2}+x_{1}x_{3}+x_{2}x_{3},\;\;x_{1}x_{2}x_{3},\\
&& x_{10}+x_{13},\;\;x_{10}x_{13},\;\;x_{6}+x_{9}+x_{11}+x_{12}+x_{14}+x_{15},\\
&& x_{11}x_{12}+x_{6}x_{14}+x_{9}x_{15},\;\;x_{9}x_{11}+x_{6}x_{12}+x_{14}x_{15},\\
&& x_{6}x_{11}+x_{9}x_{12}+x_{9}x_{14}+x_{12}x_{14}+x_{6}x_{15}+x_{11}x_{15},\\
&& x_{6}x_{9}x_{14}+x_{6}x_{11}x_{14}+x_{11}x_{12}x_{14}+x_{6}x_{9}x_{15}
%&& \qquad\mbox{}
+x_{9}x_{12}x_{15}+x_{11}x_{12}x_{15},\\
&& x_{6}^6+x_{9}^6+x_{11}^6+x_{12}^6+x_{14}^6+x_{15}^6,\;\;x_{4},\;\;x_{5}+x_{7}+x_{8},\\
&& x_{5}x_{7}+x_{5}x_{8}+x_{7}x_{8},\;\;x_{5}x_{7}x_{8}
   \end{eqnarray*}
   There are 1728 secondary invariants of maximal degree 17, among which are 76 irreducible
   secondary invariants of maximal degree 4.

\item A $18$--dimensional representation of $S_3$ is given by the matrices
   \begin{eqnarray*}
    M_1 &=& \left(e_{2}e_{1}e_{3}e_{4}e_{12}e_{10}e_{7}e_{11}e_{14}e_{6}e_{8}e_{5}e_{15}e_{9}e_{13}e_{17}e_{16}e_{18}\right)
%    \begin{smallmatrix}
%0&1&0&0&0&0&0&0&0&0&0&0&0&0&0&0&0&0\\
%1&0&0&0&0&0&0&0&0&0&0&0&0&0&0&0&0&0\\
%0&0&1&0&0&0&0&0&0&0&0&0&0&0&0&0&0&0\\
%0&0&0&1&0&0&0&0&0&0&0&0&0&0&0&0&0&0\\
%0&0&0&0&0&0&0&0&0&0&0&1&0&0&0&0&0&0\\
%0&0&0&0&0&0&0&0&0&1&0&0&0&0&0&0&0&0\\
%0&0&0&0&0&0&1&0&0&0&0&0&0&0&0&0&0&0\\
%0&0&0&0&0&0&0&0&0&0&1&0&0&0&0&0&0&0\\
%0&0&0&0&0&0&0&0&0&0&0&0&0&1&0&0&0&0\\
%0&0&0&0&0&1&0&0&0&0&0&0&0&0&0&0&0&0\\
%0&0&0&0&0&0&0&1&0&0&0&0&0&0&0&0&0&0\\
%0&0&0&0&1&0&0&0&0&0&0&0&0&0&0&0&0&0\\
%0&0&0&0&0&0&0&0&0&0&0&0&0&0&1&0&0&0\\
%0&0&0&0&0&0&0&0&1&0&0&0&0&0&0&0&0&0\\
%0&0&0&0&0&0&0&0&0&0&0&0&1&0&0&0&0&0\\
%0&0&0&0&0&0&0&0&0&0&0&0&0&0&0&0&1&0\\
%0&0&0&0&0&0&0&0&0&0&0&0&0&0&0&1&0&0\\
%0&0&0&0&0&0&0&0&0&0&0&0&0&0&0&0&0&1
%              \end{smallmatrix}\right) 
              \\
    M_2 &=& \left(e_{1}e_{3}e_{2}e_{14}e_{8}e_{7}e_{6}e_{5}e_{9}e_{10}e_{15}e_{13}e_{12}e_{4}e_{11}e_{16}e_{18}e_{17}\right)
%    \begin{smallmatrix}
%1&0&0&0&0&0&0&0&0&0&0&0&0&0&0&0&0&0\\
%0&0&1&0&0&0&0&0&0&0&0&0&0&0&0&0&0&0\\
%0&1&0&0&0&0&0&0&0&0&0&0&0&0&0&0&0&0\\
%0&0&0&0&0&0&0&0&0&0&0&0&0&1&0&0&0&0\\
%0&0&0&0&0&0&0&1&0&0&0&0&0&0&0&0&0&0\\
%0&0&0&0&0&0&1&0&0&0&0&0&0&0&0&0&0&0\\
%0&0&0&0&0&1&0&0&0&0&0&0&0&0&0&0&0&0\\
%0&0&0&0&1&0&0&0&0&0&0&0&0&0&0&0&0&0\\
%0&0&0&0&0&0&0&0&1&0&0&0&0&0&0&0&0&0\\
%0&0&0&0&0&0&0&0&0&1&0&0&0&0&0&0&0&0\\
%0&0&0&0&0&0&0&0&0&0&0&0&0&0&1&0&0&0\\
%0&0&0&0&0&0&0&0&0&0&0&0&1&0&0&0&0&0\\
%0&0&0&0&0&0&0&0&0&0&0&1&0&0&0&0&0&0\\
%0&0&0&1&0&0&0&0&0&0&0&0&0&0&0&0&0&0\\
%0&0&0&0&0&0&0&0&0&0&1&0&0&0&0&0&0&0\\
%0&0&0&0&0&0&0&0&0&0&0&0&0&0&0&1&0&0\\
%0&0&0&0&0&0&0&0&0&0&0&0&0&0&0&0&0&1\\
%0&0&0&0&0&0&0&0&0&0&0&0&0&0&0&0&1&0
%              \end{smallmatrix}\right) 
   \end{eqnarray*}
   Our primary invariants are
   \begin{eqnarray*}
&&x_{1}+x_{2}+x_{3},\;\;
x_{1}x_{2}+x_{1}x_{3}+x_{2}x_{3},\;\;
x_{1}x_{2}x_{3},\\
&&x_{4}+x_{9}+x_{14},\;\;
x_{4}x_{9}+x_{4}x_{14}+x_{9}x_{14},\;\;
x_{4}x_{9}x_{14},\\
&&x_{16}+x_{17}+x_{18},\;\;
x_{16}x_{17}+x_{16}x_{18}+x_{17}x_{18},\;\;
x_{16}x_{17}x_{18},\\
&&x_{6}+x_{7}+x_{10},\;\;
x_{6}x_{7}+x_{6}x_{10}+x_{7}x_{10},\\
&&x_{6}x_{7}x_{10},\;\;
x_{5}+x_{8}+x_{11}+x_{12}+x_{13}+x_{15},\\
&&x_{5}x_{12}+x_{8}x_{13}+x_{11}x_{15},\;\;
x_{8}x_{11}+x_{12}x_{13}+x_{5}x_{15},\\
&&x_{5}x_{11}+x_{8}x_{12}+x_{5}x_{13}+x_{11}x_{13}+x_{8}x_{15}+x_{12}x_{15},\\
&&x_{5}x_{8}x_{12}+x_{5}x_{11}x_{12}+x_{5}x_{8}x_{13}+x_{11}x_{12}x_{15}+x_{8}x_{13}x_{15}+x_{11}x_{13}x_{15},\\
&&x_{5}^6+x_{8}^6+x_{11}^6+x_{12}^6+x_{13}^6+x_{15}^6
 \end{eqnarray*}
   There are 31104 secondary invariants of maximal degree 22, among which are 137 irreducible
   secondary invariants of maximal degree 4.

\item A $10$--dimensional representation of $S_5$ is given by the matrices
   \begin{eqnarray*}
    M_1 &=& \left(\begin{smallmatrix}
1 &  0  & 0  & 0  & 0  & 0  & 0 &  0 &  0 &  0\\
0 &  1 &\frac 13 &\frac 13 &\frac 13&   0 &  0&   0&   0&   0\\
0 &  0 &\frac 13&-\frac 23&-\frac 23 &  0&   0 &  0&   0&   0\\
0 &  0&-\frac 23& \frac 13&-\frac 23&   0&   0&   0&   0&   0\\
0 &  0&-\frac 23&-\frac 23& \frac 13&   0&   0&   0&   0&   0\\
0 &  0&   0&   0&   0&   1&   0&   0&   0&   0\\
0 &  0&   0&   0 &  0 &  0 &  0 &  0&   1 &  0\\
0 &  0&   0 &  0 &  0 &  0 &  0 &  0 &  0 &  1\\
0 &  0 &  0 &  0 &  0 &  0 &  1 &  0 &  0 &  0\\
0 &  0 &  0 &  0 &  0 &  0  & 0 &  1 &  0 &  0
              \end{smallmatrix}\right) \\
    M_2 &=& \left(\begin{smallmatrix}
1 &  0 &  0 &  0 &  0 &  0 &  0&   0 &  0  & 0\\
0 &  0 &\frac 13&-\frac 23&-\frac 23&   0&   0&  0 &  0 &  0\\
0 &  0&-\frac 23& \frac 13&-\frac 23&   0 &  0&  0 &  0 &  0\\
0 &  0&-\frac 23&-\frac 23& \frac 13&   0 &  0&   0 &  0 &  0\\
0 &  1& \frac 13& \frac 13& \frac 13&   0 &  0&   0 &  0 &  0\\
0 &  0 &  0 &  0 &  0&  -1&  -1&   1&   1 &  0\\
0 &  0 &  0 &  0 &  0 & -1 &  0&   0&   0 &  1\\
0 &  0 &  0 &  0&   0 & -1&   0 &  1 &  0  & 0\\
0 &  0 &  0 &  0 &  0 &  0&  -1&  0  & 0 &  0\\
0  & 0 &  0 &  0 &  0 &  0&  -1&  1 &  0 &  0
              \end{smallmatrix}\right)
   \end{eqnarray*}
   We are not listing the primary invariants here, as they are too big polynomials.
   There are 720 secondary invariants of maximal degree 22, among which are 46 irreducible
   secondary invariants of maximal degree 9.
\end{enumerate}

Examples~(2), (3) and (9) belong to a very interesting class of examples that was shown to us 
by G. Kemper~\cite{KemperPersonal}.
For $n\in \mathbb N$, Let $M_n$ be the set of two-element subsets of $\{1,...,n\}$. Then,
one studies the obvious $S_n$ action on $M_n$ (or similarly, the obvious $A_n$ action), 
and one can try to compute
the invariant ring $\mathbb Q[M_n]^{S_n}$ (resp. $\mathbb Q[M_n]^{A_n}$). 

The $10$--dimensional representation of $S_5$ in Example~(9) is a surprisingly challenging
problem. To simplify the computations, Kemper provided a decomposition of the representation
into a direct sum of a $1$--, a $4$-- and a $5$--dimensional representation. 
Without ad-hoc methods, the computation of secondary invariants for that problem has been beyond reach. 
The procedure \texttt{(Ir\-re\-du\-cible)Se\-con\-dary\-Inva\-ri\-ants} of \textsc{Magma}~V2.13-8 breaks immediately, 
since it requests 55.62~GB memory, while the memory limit of our computer is 16~GB. Our algorithm 
\texttt{ir\-red\underline{\ }se\-con\-da\-ry\underline{\ }char0} 
in \textsc{Singular} version 3-0-2 exceeds the limit of 16~GB while computing secondary 
invariants in degree $8$. 

The total number of secondary invariants in Example~(9) is not particularly large.
The difficulties in Example~(9) come from the fact that there are \emph{irreducible} secondary invariants of rather high degrees. 

\subsection{Comparison}

We describe here how different algorithms perform on Examples~(1) up to~(9).
All computations had been done on a Linux x86\underline{\ }64 platform with two AMD Opteron 
248 processors (2,2 GHz) and a memory limit of 16~GB.
The computation of primary invariants is not part of our tests. Hence, in each example we use the
same primary invariants for all considered implementations.
We compare the following implementations:
\begin{enumerate}
\item \texttt{secondary\underline{\ }char0} as in \textsc{Singular} release 2-0-6. 
 In Table~\ref{tab:bench}, we refer to it as \lq\lq\textsc{Singular} (1998)\rq\rq.
\item \texttt{secondary\underline{\ }char0} as in \textsc{Singular} release 3-0-2, whith a small refinement. 
 In Table~\ref{tab:bench}, we refer to it as \lq\lq\textsc{Singular} (all sec.)\rq\rq.
\item \texttt{irred\underline{\ }secondary\underline{\ }char0}, as in
 \textsc{Singular} release 3-0-2, with a small refinement. 
 In Table~\ref{tab:bench}, we refer to it as \lq\lq\textsc{Singular} (irr. sec.)\rq\rq.
\item \texttt{SecondaryInvariants} in \textsc{Magma}~V2.13-8. 
\end{enumerate}

Implementation~(1) is due to A. Heydtmann~\cite{Heydtmann} (1998)
and has been part of \textsc{Singular} up to release~3-0-1. 

Implementations~(2) and (3) are our implementations
of the \textsc{Improved New Algorithm} explained in Section~\ref{sec:irred_secondary}.
They are part of \textsc{Singular}~3-0-2, released in Juli, 2006.
Here, we test a slightly improved version, that
saves memory when generating irreducible secondary invariants in high degrees. 
However, this only affects example~(9); the performance in the other eight examples 
remains essentially the same,
as the degrees of their irreducible secondary invariants are not high enough.

\begin{table}
\begin{tabular}[t]{l|c|c|c|c|}
Algorithm:& (1) & (2) & (3) & (4) \\
 &\textsc{Singular}  &\textsc{Singular}&\textsc{Singular} & \textsc{Magma} \\
&(1998) &(all sec.) & (irr. sec.) &                \\
 \hline
Expl. (1)&0.55 s   & 0.05 s & 0.03 s  &0.05 s     \\
         &8.62 MB  & 1.49 MB& 1.0 MB  &10.3 MB    \\
\hline
Expl. (2)&0.05 s   & 0.04 s & 0.04 s  &0.01 s     \\
         &0.99 MB  & 0.96 MB& 0.97 MB &7.05 MB    \\
\hline
Expl. (3)&0.48 s   & 0.33 s & 0.3  s  &0.19 s      \\
         &2.97 MB  & 1.95 MB& 1.96 MB&8.96 MB     \\
\hline
Expl. (4)&6.55 s   & 0.63 s & 0.32 s  &0.48 s     \\
         &12.29 MB & 2.47 MB& 2.97 MB &9.09 MB    \\
\hline
Expl. (5)&18.15 s  &10.53s  &9.69 s   &6.66 s     \\
         &45.79 MB &10.61 MB&17.0 MB  &31.82 MB   \\
\hline
Expl. (6)&$> 984$ m&100.4 s &16.55 s  &118.51 s   \\
         &$>167$ MB&110.0 MB&39.0 MB  &54.0 MB    \\
\hline
Expl. (7)&---      &268.9s  &20.94 s  &$> 7$ h    \\
         &---      &872.7 MB&35.1 MB  &$> 15$ GB  \\
\hline
Expl. (8)&---      &$>10$ h & 50.7 m  & ---       \\
         &---      &$>10$ GB& 3.36 GB &(259.5 GB) \\
\hline
Expl. (9)&---      & 6.42 h & 99.2 m  & ---       \\
         &---      &10.74 GB& 7.35 GB & (55.62 GB) \\ 
\hline
\end{tabular}
\medskip
\caption{Comparison of different implementations}
\label{tab:bench}
\end{table}

Implementation~(4) is due to A. Steel, based on~\cite{KemperSteel} or
\cite{Kemper} or~\cite{DerksenKemper}. We consider here the \textsc{Magma}-version V2.13-8, released in October, 2006.
There is also a function \texttt{Ir\-re\-du\-cible\-Se\-con\-da\-ry\-In\-va\-riants} in \textsc{Magma}, but
computation time and memory consumption are essentially the same, in our examples. So, for the sake of simplicity,
we do not provide separate timings for that function.

Note that, after posting the first version of this manuscript, there was a new release
of \textsc{Magma} containing an algorithm that G. Kemper developed in 2006. However, it
seems that Kemper did not describe his algorithm in a paper yet.
Meanwhile we implemented another, completely different algorithm in \textsc{Singular}. It will be part of \textsc{Singular} release 3-0-3 and often works much faster. E.g., it
can compute Example~(8) in $1.06$ seconds. We describe this algorithm 
in~\cite{KingInvariantAlgebra} and also provide there comparative benchmarks using the 
new versions of \textsc{Singular} and \textsc{Magma}.

Interestingly, in contrast to the corresponding \textsc{Magma} functions, 
\texttt{ir\-red\underline{\ }se\-con\-da\-ry\underline{\ }char0} often works much faster and needs much less memory 
than \texttt{se\-con\-da\-ry\underline{\ }char0}; see Examples (1) and (6)--(9). 
However, this is not always the case, as can be seen in Examples~(4) and~(5).

In Table~\ref{tab:bench}, \lq\lq ---\rq\rq\ means that the computation fails since
the process exceeds the memory limit; in examples (8) and (9), \textsc{Magma} requests
the amount of memory that we indicate in round brackets.
In some cases, we stopped the computation when it was clear that it takes too
much time; this is indicated in the table by \lq\lq$>...$\rq\rq.

In conclusion, our benchmarks provide some evidence that the \textsc{Improved New Algorithm} has great advantages in the computation of invariant rings with many secondary invariants. Here, it marks a dramatic improvement compared with previous algorithms in \textsc{Singular} or algorithms in \textsc{Magma}. 
In 3 of our 9 examples,
it is the only algorithm that terminates in reasonable time with a memory limit of 16 GB. A particular benefit or our algorithm is that the computation of irreducible secondary invariants does not involve the explicit computation of \emph{reducible} secondary invariants, which may save resources.

\subsection*{Acknowledgement}
I'm grateful to Gregor Kemper for providing me with the data of Example~(9).
I owe thanks to Gregor Kemper and Nicolas Thi\'ery for their comments on this 
manuscript.

%%%%%%%%%%%%%%%%%%%%%%%%%%%%%%%%%%%%%%%%%%%%%%%%%%%%%%%%%%%%%%%%%%%%%%%%%%%%%

\end{document}